\documentclass[12pt]{article}
\usepackage{graphicx}
\font\fmtitle=cmbx12 scaled \magstep2
\font\text=cmr10 at 12 truept
 at 11 truept
 at 10 truept

\begin{document}
 \thispagestyle{empty}
 \markboth
 {\it \centerline {Proceedings MATHAI-75}}
 {\it \centerline{Rudolf GORENFLO and Francesco MAINARDI}}  
\def\pni{\par \noindent}
\def\vsh{\vskip 0.25truecm\noindent}
\def\vs{\vskip 0.5truecm}
\def\vvs{\vskip 1.0truecm}
\def\vvvs{\vskip 1.5truecm}
\def\vsp{\vsh\pni}
\def\vsn{\vsh\pni}
\def\cen{\centerline}
\def\ra{\item{(a)\ }} \def\rb{\item{(b)\ }}   \def\rc{\item{(c)\ }}
\def\eg{{e.g.}\ } \def\ie{{i.e.}\ }
\def\sg{\hbox{sign}\,}
\def\sgn{\hbox{sign}\,}
\def\sign{\hbox{sign}\,}
\def\e{\hbox{e}}
\def\exp{\hbox{exp}}
\def\ds{\displaystyle}
\def\dis{\displaystyle}
\def\q{\quad}    \def\qq{\qquad}
\def\lan{\langle}\def\ran{\rangle}
\def\l{\left} \def\r{\right}
\def\lra{\Longleftrightarrow}
\def\arg{{\rm arg}}
\def\argz{{\rm arg}\, z}
\def\argG{{x^2/ (4\,a\, t)}}
\def\d{\partial}
 \def\dr{\partial r}  \def\dt{\partial t}
\def\dx{\partial x}   \def\dy{\partial y}  \def\dz{\partial z}
\def\rec#1{\frac{1}{#1}}
\def\log{{\rm log}\,}
\def\erf{{\rm erf}\,}     \def\erfc{{\rm erfc}\,}
\def\upphi{\phi}
\def\upkappa{\kappa}
\def\updelta{\delta}
\def\upalpha{\alpha}
\def\upbeta{\beta}
\def\uplambda{\lambda}
\def\NN{{\rm I\hskip-2pt N}}
\def\MM{{\rm I\hskip-2pt M}}
\def\RR{\vbox {\hbox to 8.9pt {I\hskip-2.1pt R\hfil}}\;}
\def\CC{{\rm C\hskip-4.8pt \vrule height 6pt width 12000sp\hskip 5pt}}
\def\II{{\rm I\hskip-2pt I}}
\def\erf{{\rm erf}\,}   \def\erfc{{\rm erfc}\,}
\def\exp{{\rm exp}\,} \def\e{{\rm e}}
\def\ss{{s}^{1/2}}   
\def\N{\bar N}  
\def\ss{{s}^{1/2}} 
\def\stt{{\sqrt t}}
\def\lst{{\lambda \,\stt}}
\def\Et{{E_{1/2}(\lst)}}
\def\u{\widetilde{u}}
\def\ul{\widetilde{u}} 
\def\uf{\widehat{u}} 
\def\bar{\widetilde}
\def\A{{\mathcal {A}}}
\def\L{{\cal L}} 
\def\F{{\cal F}} 
\def\M{{\cal M}}  
\def\P{{\cal P}}  
\def\Fdiv{\,\stackrel{{\cal F}} {\leftrightarrow}\,}
  \def\Ldiv{\,\stackrel{{\cal L}} {\leftrightarrow}\,}
  \def\Mdiv{\,\stackrel{{\cal M}} {\leftrightarrow}\,}
\def\args{(x/ \sqrt{a})\, s^{1/2}}
\def\argsa{(x/ \sqrt{a})\, s^{\beta}}
\def\barr{\widetilde}
\def\argG{ x^2/ (4\,a\, t)}
\def\G{{\cal {G}}}
\def\Gc{{\cal {G}}_c}	\def\Gcs{\barr{\Gc}} 
\def\Gs{{\cal {G}}_s}	\def\Gss{\barr{\Gs}} 
\def\f{\bar{f}}
\def\g{\bar{g}}
\def\u{\bar{u}}
\def\be{\begin{equation}}
\def\ee{\end{equation}}

\centerline{{\fmtitle On the fractional Poisson process}}
\vskip 0.20truecm  
\centerline{{\fmtitle and the discretized stable subordinator }}
\vskip 0.20truecm
\begin{center}
{{\bf Rudolf GORENFLO$^1$ and Francesco MAINARDI$^2$}}
 \vskip 0.25truecm
$^1$ Dept. of Mathematics \& Informatics, Free University Berlin, Germany
\\ E-mail: {\tt gorenflo@mi.fu-berlin.de}    
\vskip 0.25truecm
$^2$ Dept. of Physics \& Astronomy, University of Bologna, and INFN,  Italy
\\ E-mail: {\tt francesco.mainardi@unibo.it} \ {\tt francesco.mainardi@bo.infn.it}    
\vskip 0.25truecm
 \end{center}
\noindent
{\bf Paper dedicated to Professor A.M.Mathai on the occasion of his 
80-th anniversary,
 published in Axioms, Vol 4, pp. 321--344 (2015). 
DOI:10.3390/axioms4030321} 
\vskip 0.25 truecm
\centerline{{\bf Abstract}}
\noindent
We consider the {\it renewal counting number process} $N=N(t)$ as a forward march over the non-negative 
integers with independent identically distributed waiting times. 
We embed the values of the counting numbers $N$ in a "pseudo-spatial" non-negative half-line $x \ge 0$ 
and observe that for  physical time likewise we have $t \ge 0$. Thus  we apply the Laplace transform 
with respect to both variables  $x$ and  $t$.
Applying then a modification of the Montroll-Weiss-Cox formalism of continuous time random walk 
we obtain the essential characteristics of a renewal process in the transform domain and, 
if we are lucky, also in the physical domain. 
The process  $t=t(N)$ of accumulation of waiting times is inverse to the  counting number process, 
in honour of the Danish mathematician and telecommunication engineer A.K. Erlang we call it the {\it Erlang  process}. 
It yields the probability of exactly $n$ renewal events in the interval $(0, t]$. 
We apply our Laplace-Laplace formalism to the fractional  Poisson process whose waiting times are 
of Mittag-Leffler type and to a renewal process whose waiting tímes are of Wright type. 
The process of Mittag-Leffler type includes as a limiting case the classical Poisson process,   the process of
Wright type represents the discretized stable subordinator and a re-scaled version of it was used 
in our method of parametric subordination of time-space fractional diffusion processes. 
Properly rescaling  the counting number process $N(t)$ and the Erlang process $t(N)$ 
 yields as diffusion limits the inverse stable and the  stable subordinator, respectively.  
\newpage 
\centerline{{\bf Contents}} 
\vskip 0.25truecm \noindent
\\ 1. Introduction
\\ 2. Elements of renewal theory and CTRW
\\ 3. The Poisson process and its fractional generalization
\\ 4. The stable subordinator and  the Wright process
\\ 5. The diffusion limits for the fractional Poisson and the Wright processes
\\ 6. Conclusions 
\\  \phantom{7.} Acknowledgments
 \\ \phantom{8.} Appendix A: Operators, transforms and special functions
 \\ \phantom{8.} Appendix B: Collection of results
\\ \phantom{9.} References

\section{Introduction}

Serious studies of the fractional generalization of the Poisson process - 
replacement of the exponential waiting time distribution by a distribution given via a Mittag-Leffler  
function with modified argument - have been started around the turn of the millenium, 
and since then many papers on its various aspects have appeared. 
There are in the literature many papers on this generalization 
where the authors have outlined a number of aspects and definitions, see e.g.
 Repin and Saichev  (2000) \cite{Repin-Saichev_RQR00}, 
Wang et al.  (2003,2006) \cite{Wang-Wen_CSF03,Wang-Wen-Zhang_CSF06},
Laskin (2003,2009) \cite{Laskin_CNSNS03,Laskin_JMP09},
Mainardi et al. (2004)\cite{Mainardi-et-al_VIETNAM04}, 
Uchaikin et al. (2008) \cite{Uchaikin-et-al_IJBC08},  
Beghin and Orsingher  (2009) \cite{Beghin-Orsingher_EJP09}, 
Cahoy et al. (2010) \cite{Cahoy-et-al_JSPI10},
Meerschaert et al. (2011) \cite{M3_FPP},
Politi et al. (2011) \cite{Politi-Scalas_EPL11},
Kochubei  (2012) \cite{Kochubei_IEOT11},
so that it seems impossible to list them all exhaustively. 
 However, in effect this generalization was used already in 1995:  Hilfer and Anton 
 \cite{Hilfer-Anton_PRE95} (without saying it in our words) 
 showed that the Fractional Kolmogorov-Feller equation (replacement of the first order time derivative 
 by a fractional derivative of order between 0 and 1)
 requires the underlying random walk to be subordinated to a renewal process 
 with Mittag-Leffler waiting time. 
\vsp
 Here we will present our formalism for obtaining the essential characteristics of a generic renewal 
 process and apply it to get those of the fractional Poisson counting process and its inverse, 
 the fractional Erlang process. 
 Both of these comprise as limiting cases the corresponding well-known non-fractional processes 
 that are based on exponential waiting time. 
 Then we will analyze an alternative renewal process, that we call the "Wright process",
 investigated by Mainardi et al (2000), (2005), (2007) 
 \cite{Mainardi-Raberto-Gorenflo-Scalas_PhysicaA00,
 Mainardi-Gorenflo-Vivoli_FCAA05,Mainardi-Gorenflo-Vivoli_JCAM07},
  a process arising by discretization of the stable subordinator. 
  In it the so-called $M$-Wright function plays the essential role. 
  A scaled version of this process has been used by Barkai (2002) \cite{Barkai_ChemPhys02} 
  for approximating the time-fractional diffusion process directly by a random walk 
  subordinated to it (executing this scaled version in natural time), 
  and he has found rather poor convergence in refinement. 
  In Gorenflo et al. (2007) \cite{GorMaiViv_CSF07} we have modified the way of using 
   this discretized stable subordinator.   
  By appropriate discretization of the relevant  spatial stable process we have then obtained
  a simulation method    equivalent
to the solution of a pair of Langevin equations, 
see Fogedby (1994) \cite{Fogedby_PRE94} and Kleinhans and Friedrich(2007) \cite{Kleinhans-Friedrich_PRE07}. 
For simulation of space-time fractional diffusion one so obtains a sequence 
of precise snapshots of a true particle trajectory,  see for details 
Gorenflo et al. (2007) \cite{GorMaiViv_CSF07}, 
and also Gorenflo and Mainardi (2011, 2012) \cite{Gorenflo-Mainardi_EPJ-ST11,Gorenflo-Mainardi_METZLER11}.
\vsp
  However, we should note that already in the Sixties of the past century, 
  Gnedenko and Kovalenko  (1968) \cite{Gnedenko-Kovalenko_QUEUEING68}
  obtained in disguised form the fractional Poisson process by 
  properly rescaled infinite thinning
(rarefaction) of a renewal process  with power law waiting time.
By "disguised" we mean that they found the Laplace transform of the Mittag-Leffler
 waiting time density, but being ignorant of the Mittag-Leffler function
 they only presented this Laplace transform. 
 The  same ignorance of the Mittag-Leffler function we again  meet in a 1985 paper by 
 Balakrishnan \cite{Balakrishnan_85}, 
 who exhibited the Mittag-Leffler waiting time density in Laplace 
 disguise as essential for approximating time-fractional diffusion for which he used the 
 description in form of a fractional integro-differential equation. 
   We have shown that the Mittag-Leffler waiting time density  
   in a certain sense is asymptotically universal for power law renewal
processes, see Gorenflo and Mainardi (2008) \cite{Gorenflo-Mainardi_BAD-HONNEF08},
Gorenflo (2010) \cite{Gorenflo_PALA10}.         
\vsp
The structure of our paper is as follows.
In Section 2 we discuss the elements of the general renewal theory and the CTRW concept. 
In Section 3 we introduce the Poisson process and its fractional generalization 
then, in Section 4, the so-called Wright process  related to the stable subordinator 
and its discretization. For both processes we consider the corresponding inverse processes,
the Erlang processes.   
In Section 5 we briefly discuss the diffusion limit for all the above processes.  
 Section 6 is devoted to conclusions.
We have collected in  Appendix A  notations and terminology, 
in particular the basics on operators, integral transforms  and special functions 
required for understanding our analysis. 
Finally, we provide in Appendix B  
 an overview on the essential results.    
\vsp
For related aspects of subordination we refer the readers to our papers 
\cite{Gorenflo-Mainardi_BAD-HONNEF08,Gorenflo-Mainardi_EPJ-ST11,%
Gorenflo-Mainardi_METZLER11,Gorenflo-Mainardi_KILBAS12,GorMaiViv_CSF07}
and to papers by Bazhlekova \cite{Bazhlekova_Mathematics15},
 Meerschaert, Nane and Vellaisamy \cite{M3_FPP},
 Umarov \cite{Umarov_FCAA15}.

\section{Elements of renewal theory and CTRW}\label{sec:2}
 For the reader's convenience let us here present a brief introduction to renewal theory including the basics of continuous time random walk (CTRW).
 
 \paragraph{The general renewal process.}
By a {\it renewal process} we mean an infinite  sequence
$0=t_0<t_1<t_2<\cdots$
  of events separated by i.i.d. (independent and identically distributed) random waiting times
  $T_j=t_j-t_{j-1}$, whose probability density $\phi(t)$  is given as a function
   or generalized function in the sense of Gel'fand and Shilov \cite{Gelfand-Shilov_BOOK64}
   (interpretable as a measure) with support on the positive real axis
   $t\ge 0$,  non-negative: $\phi(t)\ge 0$, and normalized: ${\ds \int _0^\infty \!\!\phi(t)\, dt} =1$,
   but not having a delta peak at the origin $t=0$.
   The instant $t_0=0$ is not counted as an event.
   An important global characteristic of a renewal process is its mean waiting time
   $\left< T \right>  = {\ds \int _0^\infty \!\! t\,\phi(t)\, dt} $. 
   It may be finite or infinite.
   In any renewal process we can distinguish two processes, namely the {\it counting number process}
   and the process inverse to it, that we call the {\it Erlang process}.
   The instants $t_1,  t_2, t_3, \dots$ are often  called  {\it renewals}. In fact renewal theory is relevant
   in practice  of maintenance or required exchange of failed parts, e.g., light bulbs.    
\paragraph{The counting number process and its inverse.}
We are interested in the {\it counting number process} $x=N= N(t)$ 
 $$ N(t):=
   \hbox{max} \left\{n | t_n \le t \right\} = n
   \q \hbox{for}\q t_n\le t<t_{n+1}\,,
    \q n=0,1,2, \cdots,
	\eqno(2.1)$$
   where in particular $N(0)=0$.
   We ask for the counting number probabilities in $n$, evolving in $t$,
$$ p_n(t):= \P [N(t)=n]\,, \; n=0,1,2, \cdots\,.
\eqno(2.2)$$
We denote by $p(x,t)$ the sojourn density for the counting number
  having the value $x$.
   For this process  the expectation is
$$
 m(t):= \left< N(t)\right>  = 
 \sum_{n=0}^\infty n\, p_n(t)=
 \int _0^\infty \!\! x\,p(x,t)\, dx\,,
 \eqno(2.3)$$
 [since $ p(x,t)= {\ds \sum_{n=0}^\infty p_n(t)\, \delta(x-n)}$, see (2.12)]
 It provides  the mean number of  events in the half-open interval   $(0,t]$,
and is called the {\it renewal function}, see e.g. \cite{Ross_BOOK96}.
We also will look at the process $t=t(N)$, the inverse to the process $N=N(t)$,
that  we call the {\it Erlang process}
in honour of the Danish telecommunication engineer A.K. Erlang (1878-1929), 
see Brockmeyer et al. (1948) \cite{Erlang-Life_1948}.
It gives the value of time $t=t_N$ of the $N$-th renewal.
We ask for the Erlang probability densities 
$$q_n(t) = q(t, n)\,, \; n=0,1,2,\ldots \eqno(2.4)  $$
For every $n$ the function $q_n(t) = q(t, n)$ is a density in the variable 
of time having value $t$ in the instant of the $n$-th event. 
Clearly, this event occurs after $n$ (original) waiting times have passed, so that 
$$ q_n(t) = \phi^{*n}(t)  \q \hbox{with Laplace transform} \q \widetilde q_n(s) = (\widetilde \phi(s)^n)\,.
\eqno(2.5)$$
In other words  the function $q_n(t) = q(t, n)$ is a probability density in the variable 
$t \ge 0$ evolving in the variable $x= n = {0, 1, 2, ...}$.
\vsp
\paragraph{The continuous time random walk.}
A {\it continuous time random walk}  (CTRW) is given by an infinite sequence of 
spatial positions
$0=x_0, x_1, x_2, \cdots$, separated by (i.i.d.) random jumps    $X_j=x_j-x_{j-1}$,
whose probability density function $w(x)$
is given as a non-negative function or generalized function (interpretable as a measure) 
with support on the real axis
$-\infty <x < +\infty$  and normalized: ${\ds \int _0^\infty \!\! w(x)\, dx} =1$,
 this random walk being subordinated to a renewal process so that we have a random process
 $x=x(t)$    on the real axis with the property $x(t)=x_n$   for $t_n\le t <t_{n+1}$, $n=0,1,2, \cdots$.
  \vsp
     We ask for the {\it sojourn probability density} $u(x,t)$
	  of a  particle wandering according to the random process $x=x(t)$
	       being in point $x$ at instant $t$.
\vsp
Let us define the following cumulative probabilities related to the probability density function $\phi(t)$
$$
	\Phi(t)  = \int_0^{t+} \!\! \phi (t')\, dt'\,, \quad 	
	\Psi(t)  = \int_{t+}^\infty \!\! \phi (t')\, dt' = 1-\Phi(t)\,.
\eqno(2.6)$$
For definiteness, we take $\Phi(t)$ as right-continuous, $\Psi(t)$ as left-continuous.
 When the non-negative random variable  represents
 the lifetime of a technical system, it is common to
 call $\Phi(t):= \P \left(T \le  t\right) $  the {\it failure probability}
 and $\Psi(t) := \P \left(T > t\right)$
  the {\it survival probability}, because $\Phi(t)$ and $\Psi(t)$ are
the respective probabilities that the system does or does not fail
in $(0, t]$. These terms, however,  are commonly adopted for  any
renewal process.
\vsp
 In the Fourier-Laplace domain we have
 $$\widetilde \Psi(s)= \frac {1-\widetilde \phi(s)}{s}\,,
 \eqno(2.7) $$ 
  and the famous Montroll-Weiss solution formula for a CTRW, see 
  \cite{Montroll-Weiss_JMP65,Weiss_BOOK94}
$$
\widehat{\widetilde{u}} (\kappa,s)
={\ds \frac{1-\widetilde{\phi}(s)}{s} \,\sum_{n=0}^\infty  \left(\widetilde \phi(s)\,\widehat w(\kappa) \right)^n}
={\ds \frac{1-\widetilde{\phi}(s)}{s} \, \frac{1}{1- \widetilde{\phi}(s)\,\widehat{w}(\kappa)}}
\,.
\eqno (2.8)$$  
In our special situation   the jump density has support only on the positive semi-axis   $x\ge 0$
and thus,  by replacing 
 the Fourier transform by the Laplace transform we obtain the Laplace-Laplace solution
$$  
\widetilde{\widetilde{u}} (\kappa,s)
 ={\ds \frac{1-\widetilde{\phi}(s)}{s}\,\sum_{n=0}^\infty  \left(\widetilde\phi (s)\,\widetilde w(\kappa) \right)^n}
 {\ds =\frac{1-\widetilde{\phi}(s)}{s} \, \frac{1}{1- \widetilde{\phi}(s)\,\widetilde{w}(\kappa)}}
\,.
\eqno(2.9) $$ 
Recalling from Appendix  the definition of convolutions,
 in the physical domain
  we have for the solution $u(x,t)$  the   {\it Cox-Weiss series}, see \cite{Cox_RENEWAL67,Weiss_BOOK94},
$$ u(x,t) = \left(\Psi \,*\, \sum_{n=0}^\infty  \phi^{*n}\, w^{*n}  \right)(x,t)\,.
\eqno(2.10)$$
This formula has an intuitive meaning: Up to and including instant $t$,  there have occurred
0 jumps, or 1 jump, or 2 jumps, or $\dots$, and if the    last jump has occurred at instant $t'<t$,
the wanderer is resting there for a duration $t-t'$. 
	   \vsp
From the rich literature on the concept of CTRW and its applications we recommend to study the surveys by Metzler and Klafter 
\cite{Metzler-Klafter_PhysRep00,Metzler-Klafter_JPA04} and 
the original article 
by Chechkin, Hofmann and Sokolov \cite{Chechkin-et-al_PRE09}. 
\vsp
\paragraph{Renewal process as a special CTRW}
The essential trick of what follows consists  in a rather non-conventional  use of the CTRW concept. We treat renewal processes as continuous time random walks
with waiting time density  $\phi(t)$ and special jump density $w(x)=\delta(x-1)$
corresponding to the fact that the counting number $N(t)$  increases by 1 at each positive event instant $t_n$.
We then have $\widetilde w(\kappa)= \exp(-\kappa)$  and  get for the counting number process  $N(t)$
  the sojourn density	in the transform domain ($s\ge 0$, $\kappa \ge 0$),
$$
\widetilde{\widetilde{p}} (\kappa,s)
={\ds \frac{1-\widetilde{\phi}(s)}{s} \,\sum_{n=0}^\infty  \left(\widetilde \phi (s)\right)^n\,\e^{-n\kappa}} 
={\ds \frac{1- \widetilde{\phi}(s)}{s}\,\frac{1}{1-\widetilde{\phi}(s)\, \e^{-\kappa}}}
\,.
\eqno(2.11) $$  
 From this formula we can find formulas for the renewal function   $m(t)$
and the probabilities $p_n(t)=P\{N(t)=n\}$.
Because  $N(t)$ assumes as values only the non-negative integers,
the sojourn density $p(x,t)$  vanishes if $x$ is not equal to one of these,
but has a delta peak of height $p_n(t)$  for $x=n$ ($n=0,1,2,3,\cdots$).
Hence
$$ p(x,t)= \sum_{n=0}^\infty p_n(t)\, \delta(x-n)\,.
\eqno(2.12) $$ 
 Inverting   (2.11)  with respect to $\kappa$ and $s$ as
 $$  p(x,t)= \sum_{n=0}^\infty \left(\Psi \,*\, \phi^{*n}\right)(t)\, \delta(x-n)\,,
  \eqno(2.13) $$  
   we identify
$$ p_n(t)= \left( \Psi\,*\, \phi^{*n}\right) (t)\,.
\eqno (2.14)$$   
According to the theory of Laplace transform we conclude from Eqs. (2.2) and (2.12)
$$%
m(t)= {\ds - \frac{\partial}{\partial \kappa} \left.\widetilde p(\kappa,t)\right|_{\kappa=0}}
= {\ds \left.\left(\sum_{n=0}^\infty n\,p_n(t)\, \e^{-n\kappa}\right)\right|_{\kappa=0}}
 = {\ds \sum_{n=0}^\infty n\, p_n(t)}\,.
\eqno(2.15) $$ 
a result naturally expected, and
$$
\widetilde m(s)
= {\ds \sum_{n=0}^\infty n\, \widetilde p_n(s)}=
{\ds \widetilde\Psi(s)\, \sum_{n=0}^\infty n \, \left(\widetilde \phi(s)\right)^n} 
= {\ds \frac{\widetilde \phi(s)}{s\left(1-\widetilde \phi(s)\right)}}\,,
\eqno(2.16)$$ 
thereby using the identity
$$ \sum_{n=0}^\infty nz^n = \frac{z}{(1-z)^2}\,, \quad |z|<1\,.$$
Thus we have found in the  Laplace domain the reciprocal pair of relationships
$$  \widetilde m(s)= \frac{\widetilde \phi(s)}{s(1-\widetilde \phi(s))}\,,
  \quad
  \widetilde \phi(s)= \frac{s\,\widetilde m(s)}{1+ s\,\widetilde m(s))}\,,
\eqno(2.17) $$ 
  saying that the waiting time density and the renewal function mutually determine each other uniquely.
  The first formula of Eq. (2.17) can also be obtained  as the value at  $\kappa=0$
   of the negative derivative for $\kappa=0$    of the last expression in Eq. (2.11).
  Eq. (2.17) implies the reciprocal pair of relationships in the physical domain
$$	
	m(t) = \int_0^t [1 + m(t-t')]\, \phi(t')\, dt'\,, \;
	m^\prime (t) =  \int_0^t [1 + m^\prime(t-t')]\, \phi(t')\, dt'\,.
\eqno(2.18)$$	
The first of these equations usually is called the {\it renewal equation}.
\vsp
Considering, formally, the counting number process $N = N(t)$ as CTRW (with jumps fixed to unit jumps 1), 
 $N$ running increasingly through the non-negative integers $x = {0, 1, 2, ...}$, 
happening in natural time $t \in [0,  \infty)$, 
we note that in the Erlang process $t = t(N)$, 
the roles of $N$ and $t$ are interchanged. The new "waiting time density" now is 
$w(x)= \delta(x-1)$, the new "jump density" is $\phi(t)$. 
\vsp
It is illuminating to consciously perceive the relationships for $t\ge 0$, $n=0,1,2, \ldots$,
between the counting number probabilities $p_n(t)$ and the Erlang densities $q_n(t)$.
For Eq. (2.5) we have $q_n(t)= \phi^{*n}(t)$, and then by (2.14)
$$ p_n(t) = \left(\Psi \,*\,q_n\right)(t) = \int_0^t \left( q_n(t')- q_{n+1}(t)\right)\, dt'\,. \eqno(2.19)$$
We can also express the $q_n$ in another way by the $p_n$. Introducing the cumulative probabilities
$ Q_n(t) = {\ds \int_0^t \!q_n(t')\,dt'}$, we have
$$Q_n(t) = \P\left( \sum_{k=1}^n T_k \le t\right)=
  \P\left(N(t)\ge  n\right)= \sum_{k=n}^\infty p_k(t)\,,\eqno(2.20)$$
  finally
 $$ q_n(t)= \frac{d}{dt} Q_n(t)=  \frac{d}{dt} \sum_{k=n}^\infty p_k(t)\,. \eqno(2.21)$$
 All this is true for $n=0$ as`well, by the empty sum convention
$  {\ds \sum_{k=1}^n T_k }=0$ for $n=0$.
 
\section{The Poisson process and its fractional \\ generalization}\label{sec:3}
  The most popular renewal process is the {\it Poisson process}.
  It is (uniquely) characterized by its {\it mean waiting time} $1/\lambda$   (equivalently by its {\it intensity} $\lambda$),
  which is a given positive number, and by its residual waiting time $\Psi(t)= \exp (-\lambda t)$
   for  $t\ge 0$, which corresponds to the  {\it waiting time density}
   $\phi(t)= \lambda \, \exp(-\lambda t)$.
	   With $\lambda=1$   we have what we call the {\it standard Poisson process}.
	  The general Poisson process arises from the standard one by rescaling the time variable $t$. 
	   \vsp
	   We generalize the standard Poisson process by replacing the
	   exponential function by  a function of Mittag-Leffler type.
	   With $t\ge 0$ and a parameter $\beta\in (0,1]$
	      we take
$$	\left\{
	\begin{array}{ll}
	\Psi(t) & \! = E_\beta(-t^\beta)\,, \\ 
	\phi(t) & \! = - {\ds \frac{d}{dt}} E_\beta(-t^\beta)
	= \beta t^{\beta-1}E^\prime_\beta(-t^\beta)  = t^{\beta-1}\, E_{\beta,\beta}(-t^\beta)
	\,.
	 \end{array}
	 \right.
\eqno(3.1)	$$  
\begin{figure}
 \includegraphics[width=.52\textwidth]{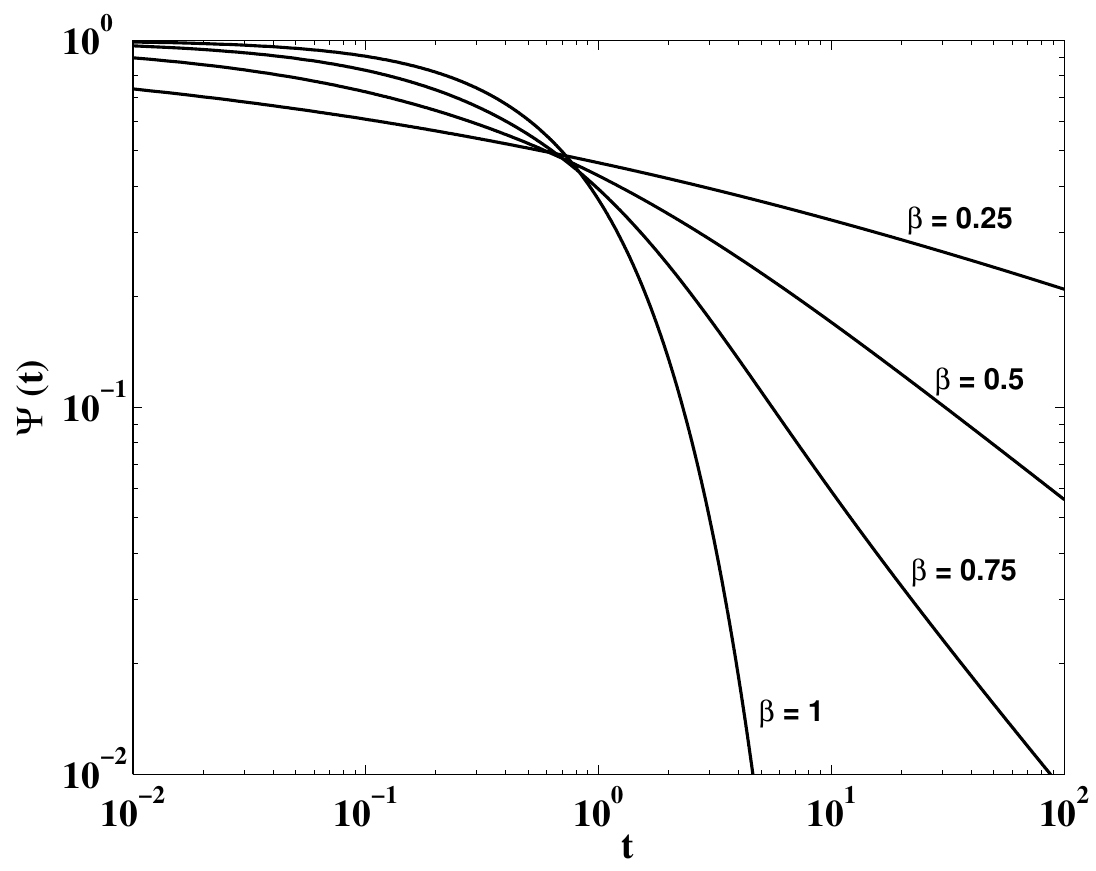}
\includegraphics[width=.52\textwidth]{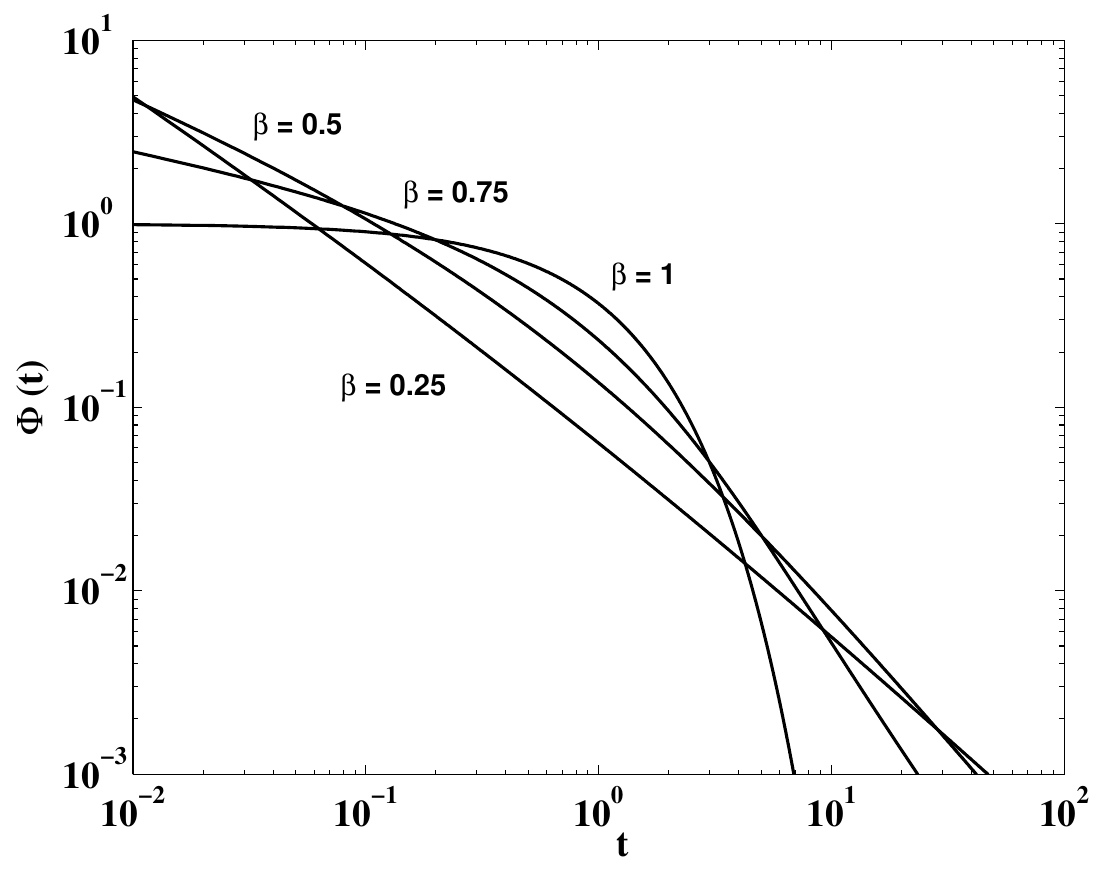}
\caption{The functions $\Psi(t)$ (left) and
$\phi(t)$ (right)
versus $t$ ($10^{-2}<t<10^{2}$) for the renewal  processes
of  Mittag-Leffler  type with $\beta = 0.25, 0.50, 0.75, 1$.}
\end{figure}
\vsp
	We call this renewal process  of Mittag-Leffler type  the {\it fractional Poisson process}, see e.g.
\cite{BDST_BOOK12,%
Beghin-Orsingher_EJP09,Cahoy-et-al_JSPI10,Gorenflo-Mainardi_KILBAS12,
Laskin_CNSNS03,Mainardi-et-al_VIETNAM04,M3_FPP,Repin-Saichev_RQR00,
Politi-Scalas_EPL11,Scalas_METZLER11},
	  and 
\cite{Uchaikin-et-al_IJBC08,Uchaikin_BOOK13},
	  or the {\it Mittag-Leffler renewal process} or the {\it Mittag-Leffler waiting time process}.
	\vsp
	To analyze it we go into the Laplace domain where we have
	$$\widetilde \Psi(s)= \frac{s^{\beta-1}}{1+ s^\beta}\,, \quad
	\widetilde \phi(s)= \frac{1}{1+ s^\beta}\,.
	\eqno(3.2)$$  
	 If there is no danger of misunderstanding we will not decorate $\Psi$ and $\phi$  with the index $\beta$.
	The special choice $\beta=1$  gives us the standard Poisson process with
	$\Psi_1(t) = \phi_1(t)= \exp(-t)$.
	\vsp
	Whereas the Poisson process has  finite mean waiting time (that of its standard version is equal to 1), the
	{\it fractional Poisson process} ($0<\beta<1$ )  does not have this property. In fact,
  $$\!	\langle T\rangle \!=\! \int_0^\infty\!\!t\,\phi(t)\, dt =
\beta \left.{\ds \frac{s^{\beta-1}}{(1+s^\beta)^2}}\right|_{s=0}
	\!=\! \left\{
	\begin{array} {ll}
	1 \,, & \beta=1\,,\\
	\infty \,, & 0<\beta<1\,.
	\end{array}
	\right.
	\eqno(3.3)$$  
Let us calculate the renewal function  $m(t)$.
Inserting $\widetilde \phi(s)= 1/(1+s^\beta)$ into Eq. (2.11) and taking
$w(x)= \delta(x-1)$ as in Section 2,
 we find for the sojourn density of the counting function $N(t)$  the expressions
$$\widetilde{\widetilde{p}} (\kappa,s)
=\frac{s^{\beta-1}}{1+s^\beta -\e^{-\kappa}}
 =  \frac{s^{\beta-1}}{1+s^\beta} \,\sum_{n=0}^\infty  \frac{\e^{-n\kappa}}{(1+s^\beta)^n}
\,,
\eqno(3.4) $$ 
and
$$\widetilde{p} (\kappa,t) = E_\beta\left(-(1-\e^{-\kappa})t^\beta\right)\,,
\eqno (3.5)$$  
and then
$$ m(t)= - \frac{\partial}{\partial \kappa} \left.\widetilde p(\kappa,t)\right|_{\kappa=0}
= \left. \e^{-\kappa} t^\beta E^\prime_\beta\left( -(1-\e^{-\kappa})t^\beta\right)\right|_{\kappa=0}\,.
\eqno(3.6)$$  
Using $E^\prime_\beta(0)= 1/\Gamma(1+\beta)$
 now yields
 $$ m(t) = \left\{
 \begin{array}{ll}
  t\,, & \beta=1\,, \\
 {\ds \frac{t^\beta} {\Gamma(1+\beta)}}\,, & 0<\beta<1\,.
 \end{array}
 \right.	
 \eqno(3.7)$$  
 This result can also be obtained by plugging $\widetilde\phi(s)= 1/(1+s^\beta)$ into
 the first equation in (2.17) 
 which yields $\widetilde m(s)=1/{s^{\beta+1}}$ and then by Laplace inversion
 Eq. (3.7).
    \vsp
Using general  Taylor expansion
$$ E_\beta(z) = \sum_{n=0}^\infty \frac{E_\beta^{(n)}}{n!} (z-b)^n\,,
\eqno(3.8)$$
in Eq. (3.5) with $b=-t^\beta$  we get
$$ \begin{array}{ll}
 \widetilde p(\kappa,t) &=
 {\ds \sum_{n=0}^\infty \frac{t^{n\beta}}{n!}\, E^{(n)}_\beta (-t^\beta)\, \e^{-n\kappa}}\,, \; \\
 p(x,t) &=
 {\ds \sum_{n=0}^\infty \frac{t^{n\beta}}{n!}\, E^{(n)}_\beta (-t^\beta)\, \delta(x-n)}\,,
 \end{array}
 \eqno(3.9)$$  
and, by comparison with Eq. (2.12), the counting number probabilities
$$		 p_n(t)= \P\{N(t)=n\} = \frac{t^{n\beta}}{n!}\, E^{(n)}_\beta (-t^\beta)\,.
\eqno(3.10)	$$	
Observing from Eq. (3.4)
$$\widetilde{\widetilde{p}} (\kappa,s)
= 	  \frac{s^{\beta-1}}{1+s^\beta} \,\sum_{n=0}^\infty  \frac{\e^{-n\kappa}}{(1+s^\beta)^n}\,,
\eqno(3.11)$$  
and inverting with respect to $\kappa$,
$$\widetilde{p} (x,s)=
\frac{s^{\beta-1}}{1+s^\beta} \,\sum_{n=0}^\infty  \frac{\delta(x-n)}{(1+s^\beta)^n}\,,
\eqno(3.12) $$  
we finally identify
$$
\widetilde p_n(s)= \frac{s^{\beta-1}}{(1+s^\beta)^{n+1}} \,\div\,
\frac{t^{n\beta}}{n!}\, E_\beta^{(n)}(-t^\beta) = p_n(t)\,.
\eqno(3.13)$$ 
En passant we have proved an often cited special case of an inversion formula 
by Podlubny (1999) \cite{Podlubny_BOOK99}, Eq. (1.80).
 \vsp
 For the  Poisson process with intensity $\lambda >0$ we have a well-known infinite system of 
 ordinary differential equations (for $t\ge 0$), see e.g. Khintchine  \cite{Khintchine_QUEUING60},
 $$ p_0(t)=\e^{-\lambda t}\,,\quad \frac{d}{dt}p_n(t)= \lambda \left(p_{n-1}(t)- p_n(t)\right)\,,
 \q n\ge 1\,,\eqno(3.14)$$  
   with initial conditions  $p_n(0)=0$, $n=1,2,\dots$,
      which sometimes even is used to define the Poisson process.
We have an analogous system of fractional differential equations for the fractional Poisson process.
In fact, from Eq. (3.13) we have
$$(1+s^\beta)\, \widetilde p_n(s) =
\frac{s^{\beta-1}}{(1+s^\beta)^n}= \widetilde p_{n-1} (s)\,.
\eqno(3.15)$$  
Hence
$$ s^\beta \, \widetilde p_n(s)= \widetilde p_{n-1}(s) - \widetilde p_{n}(s)\,,
\eqno(3.16)$$  
so in the time domain
$$ p_0(t)=E_\beta(-t^\beta)\,,\quad \,_*D_t^\beta p_n(t)= p_{n-1}(t)- p_n(t)\,,\q n\ge 1\,,
 \eqno(3.17)$$  
with initial conditions   $p_n(0)=0$, $n=1,2,\dots$, where
$\,_*D_t^\beta$ denotes the time-fractional derivative of Caputo type of order $\beta$, see Appendix A.
It is also possible to introduce and define the fractional Poisson process 
by this difference-differential system.
\vsp 
 Let us note that by solving the system  (3.17), Beghin and Orsingher in \cite{Beghin-Orsingher_EJP09}
 introduce what they call the "first form of the fractional  Poisson process"$\,$, and in \cite{M3_FPP}
Meerschaert et al. show that this process is a renewal process with 
Mittag-Leffler waiting time density as in (3.1), hence is identical with 
the fractional Poisson process.
\vsp
Up to now we have investigated the fractional Poisson counting process $N=N(t)$ and found its probabilities $p_n(t)$ in Eq. (3.10). To get the corresponding {\it Erlang probability densities}
$q_n(t)= q(t,n)$,  densities in $t$, evolving in $n=0,1,2 \ldots$, we find 
by Eq. (2.21) via telescope summation
$$ q_n(t)= \beta \frac{t^{n\beta-1}}{(n-1)!}\, E_\beta^{(n)}\left(-t^\beta\right)\,,
\quad 0<\beta \le 1\,.\eqno (3.18)$$
We leave it as an exercise to the readers to show that in Eq. (3.9) 
interchange of differentiation and summation is allowed.
\vsp
{\bf Remark} With $\beta=1$ we get the corresponding well-known results for the standard Poisson process.
The counting number probabilities are 
$$ p_n(t)= \frac{t^n}{n!}\, \e^{-t}\,,\q n=0,1,2, \ldots\,,\; t\ge 0\,, \eqno(3.19)$$
and the Erlang densities 
$$ q_n(t)= \frac{t^{n-1}}{(n-1)!}\, \e^{-t}\,,\q n=1,2,3, \ldots ,,\; t\ge 0\,.\eqno(3.20)$$
By rescalation of time we obtain  
$$ p_n(t)= \frac{(\lambda t)^n}{n!} \,\e^{-\lambda t}\,,\q n=0,1,2, \ldots,,\; t\ge 0\,, \eqno(3.21)$$
for the classical Poisson process with intensity $\lambda$
and 
$$ q_n(t)= \lambda \,\frac{(\lambda t)^{n-1}}{(n-1)!} \,\e^{-\lambda t}\,,\q n=1,2,3, \ldots,,\; t\ge 0\,.
 \eqno(3.22)$$
for the corresponding Erlang process.
\section{The stable subordinator and the Wright process}  
Let us denote by $g_\beta(t)$ the {\it extremal L\'evy stable density} of order $\beta\in (0,1]$ and support
in $t\ge 0$ whose Laplace transform is $\widetilde g_\beta(s) = \exp(-s^\beta) $, that is
$$  t\ge 0\,, \quad g_\beta(t)\,\div\,   \exp(-s^\beta)\,, \quad Re (s) \ge  0\,, \quad 0<\beta\le 1\,. \eqno(4.1)$$
The topic of L\'evy stable distributions is treated in several books on probability and stochastic processes, see e.g. 
Feller (1971) \cite{Feller_1971}, Sato (1999) \cite{Sato_99}; an overview of the analytical and graphical
aspects of the corresponding densities is found in Mainardi et al (2001) 
\cite{Mainardi-Luchko-Pagnini_FCAA01},
where an ad hoc  notation is used.
\vsp
From the Laplace transform correspondence (4.1) it is easy to derive the analytical expressions for $\beta=1/2$
(the so-called {\it L\'evy-Smirmov density}),
$ g_{1/2}(t) = {\ds{1\over 2\sqrt{\pi}}\, t^{-3/2}\,\exp(- 1/(4t))}$   
and for the limiting case $\beta=1$ (the time drift),
$g_1(t) =\delta(t-1)$, 
where $\delta$ denotes the Dirac generalized  function.
\vsp
We note that the stable  density (4.1) can be expressed in 
terms of a function of the Wright type.  
In fact, with the M-Wright function from Appendix A of this paper 
(see Appendix F of Mainardi's book \cite{Mainardi_BOOK10} for more details), we have
$$ g_\beta(t)= \frac{\beta}{t^{\beta+1}}\, M_\beta (t^{-\beta}) \,.\eqno(4.2)$$
The renewal process with waiting time density 
$$ \phi(t)= g_\beta(t) \eqno(4.3)$$
was considered in detail by Mainardi et al. (2000), (2005), (2007) 
 \cite{Mainardi-Raberto-Gorenflo-Scalas_PhysicaA00,
 Mainardi-Gorenflo-Vivoli_FCAA05,Mainardi-Gorenflo-Vivoli_JCAM07}.
 We   call this process the {\it Wright renewal process} 
 because the corresponding survival function $\Psi(t)$ and  the waiting time density $\phi(t)$
 are expressed in terms 
 of certain Wright functions. 
 So we  distinguish it from the so called {\it Mittag-Leffler renewal process}, 
 treated in the previous Section  as {\it fractional Poisson process}.
 More precisely, recalling the Wright functions from the Appendix A, we have 
 for  $t \ge 0$,
 $$   \Psi(t) = \cases{
    1-  W_{-\beta,1} \left(- {1\over t^\beta}\right),
      & $\, 0<\beta<1,$ \cr
 \Theta(t) -\Theta(t-1),
     & $\, \ \beta=1,$ \cr}
  \; \hbox{from} \;
  \widetilde \Psi(s) = {1-\e^{\, \ds -s^\beta}\over s},
\eqno(4.4) $$
$$ \phi(t) = \cases{ {1\over t}\,W_{-\beta,0} \left(- {1\over t^\beta}\right),
  & $\, 0<\beta<1,$ \cr
 \delta(t-1),
  & $\, \ \beta=1,$ \cr}
 \; \hbox{from} \;
\widetilde\phi(s)  =    \e^{\,\ds -s^\beta},
 \eqno (4.5)$$
 where $\Theta$ denotes the unit step Heaviside function.

\begin{figure}[h!]
\includegraphics[width=.52\textwidth]{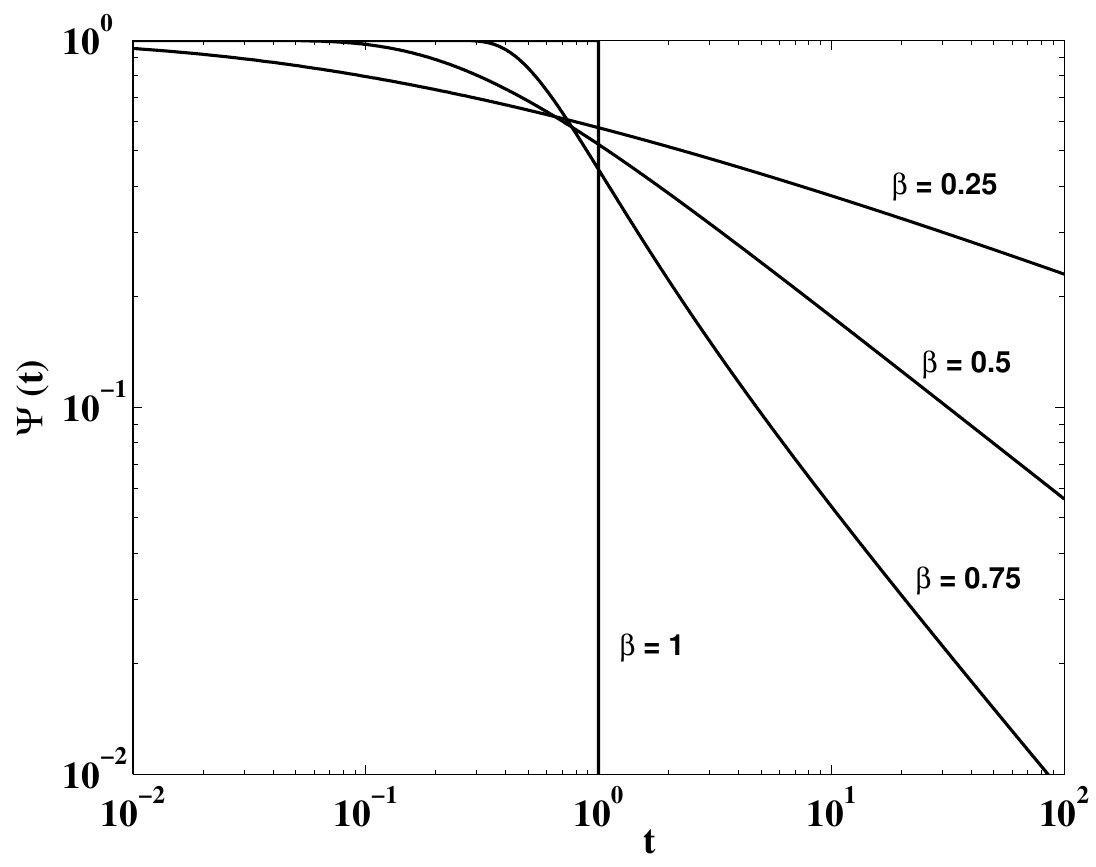}
\includegraphics[width=.52\textwidth]{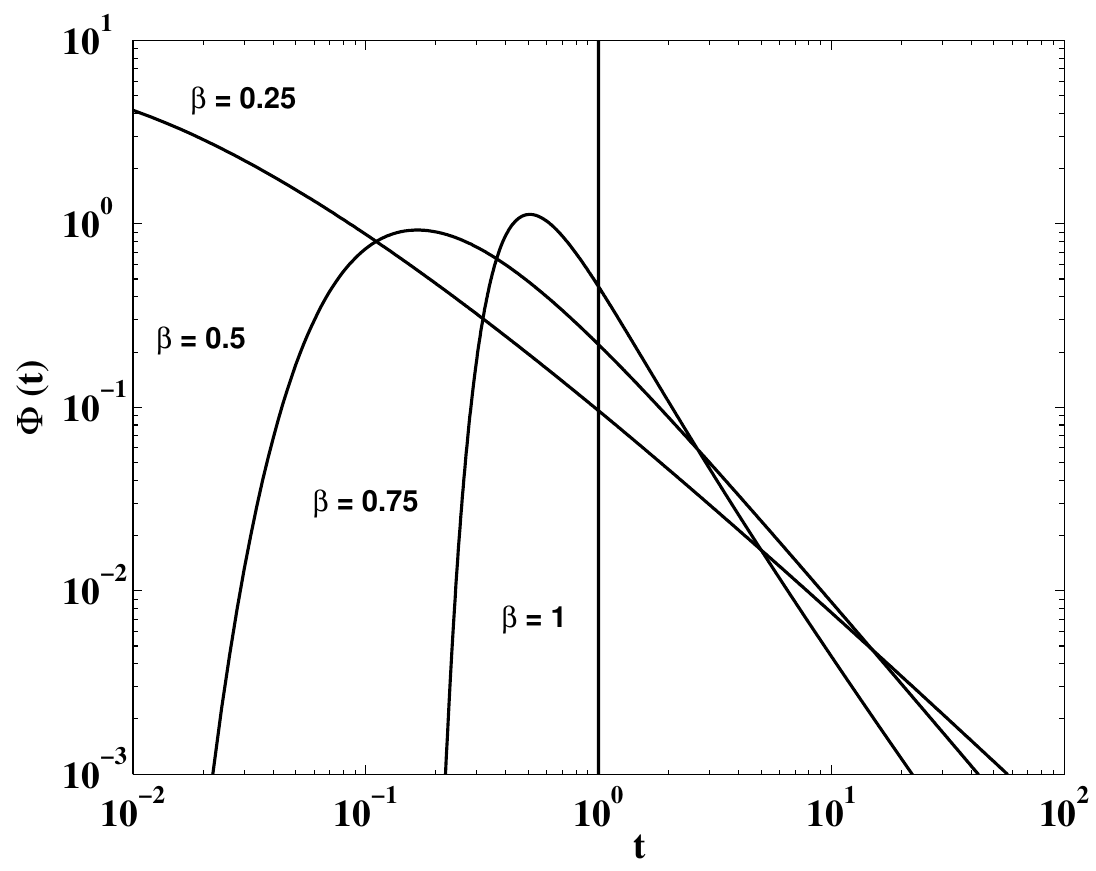}
\caption{The functions $\Psi(t)$ (left) and
$\phi(t)$ (right)
versus $t$ ($10^{-2}<t<10^{2}$) for the renewal  processes of Wright type
with $\beta = 0.25, 0.50, 0.75, 1$.
For $\beta =1$ the reader would recognize the Box function
(extended up to $t=1$) at left and the delta function
(centred in  $t=1$) at right}
\end{figure}

\vsp
 It is relevant to note the Laplace transform  connecting the two transcendental  functions $M_\beta$ and 
 $E_\beta$
 $$ M_\beta(t)\,\div\, E_\beta(-s)\,, \q 0<\beta\le 1\,. \eqno(4.6)$$   
  \vsp
By the {\it stable subordinator} of order $\beta\in (0,1]$ we mean the stochastic process 
$t=t(x)$ that has sojourn density in $t\ge 0$, evolving in $x \ge 0 $ 
provided by the Laplace transform correspondence,
$$ \!\! \widetilde f(s,x)= \e^{\ds -x s^\beta}\, \div\,  f(t,x)= x^{-1/\beta}\,g_\beta\left(x^{-1/\beta}\,t\right)
=
\frac{\beta}{t^{\beta+1}}\, x^{1+1/\beta}\, M_\beta (x t^{-\beta}) \,.
\eqno(4.7)$$ 
 This process is monotonically increasing: for this reason it is used in the context of  time change 
 and subordination in fractional diffusion processes.     
\vsp
We discretize the process  $t=t(x)$ by restricting $x$ to run through the integers 
$n = 0, 1, 2, \ldots$ . The resulting discretized version is a renewal process happening  
in pseudo-time $x \ge 0$ with jumps in pseudo-space $t\ge 0$ having density 
$g_\beta (t)$. 
Inverting this discretized stable subordinator we obtain a 
counting number process $x=N= N(t)$ with waiting time density and jump density
$$\phi(t)= g_\beta (t)\,, \quad w(x)= \delta(x-1)\,.\eqno(4.8)$$
Because here  the waiting time density is given by a function of Wright type 
we call this process the  Wright renewal process, or simply the {\it Wright process}. 
Immediately we get its {\it Erlang densities} (in $t\ge 0$, evolving in $x = n = 0, 1, 2, \ldots$)
 $$ q_n(t)= \phi^{*n}(t) \,\div \, \e^{-ns^\beta}\,, \eqno(4.9)$$
  so that, in view of (4.7) with $x=n$, 
  $$ q_n(t)=f(t,n) = n^{-1/\beta} \, g_\beta \left(n^{-1/\beta}\, t\right)\,,\eqno(4.10)$$
  In the special case $\beta=1$ we have  $q_n(t)= \delta(t-n)$.
 \vsp
 We  observe that this counting process gives us precise snapshots at $x=0, 1, 2, $ of 
 the stable subordinator 
  $t=t(x)$.
    \vsp
  Using (4.9) in  (2.14) we find the {\it counting number probabilities} in time and Laplace domain 
  $$ p_n(t) = (\Psi \,*\, \phi^{*n})(t) \,\div\, 
     \widetilde p_n(s)= {\ds \frac{1-\e^{- s^\beta} }{s}\, \e^{-n s^\beta}} = 
	 {\ds  \frac{\e^{- n s^\beta} - \e^{-(n+1) s^\beta}}{s}}\,,\eqno(4.11)$$
	  hence
	  $$p_n(t)= \int_0^t \!\left(q_n(t')-q_{n+1}(t')\right) dt'\,, \eqno(4.12)$$
	  according to (2.19).
	  \vsp
	  With the probability distribution function
	  $$ G_\beta(t) = \int_0^t g_\beta(t')\, dt'\,, \eqno(4.13)$$
	  we get 
	  $$ p_n(t)= G_\beta\left(n^{-1/\beta} t\right)-G_\beta\left((n+1)^{-1/\beta}t\right)\,.\eqno(4.14)$$
	  In the limiting case $\beta=1$ we have
	  $$ G_1(t)=\int_0^t \delta(t'-1)dt' = \left\{
	  \begin{array}{ll}
	  0 \; & \hbox{for} \; t<1\,,\\
	  1 \; & \hbox{for} \; t \ge 1\,,
	  \end{array}
	  \right.
	  \eqno(4.15)$$
	  as a function continuous from the right, and we calculate
	  $$ p_n(t)= \left\{
	  \begin{array}{ll}
	  0 \; & \hbox{for} \; 0<t<n\,, \; \hbox{and for}\; t\ge n+1\,,\\
	  1 \; & \hbox{for} \; n\le t < n+1\,.
	  \end{array}
	  \right.
	  \eqno(4.16)$$
	  For the renewal function we  obtain its Laplace transform from (2.17) 
	  $$ 
	  \widetilde m(s) =\frac{\e^{- s^\beta}} {s (1- \e^{- s^\beta})}=
	  \frac{1}{s} \sum_{n=1}^\infty \e^{- n s^\beta}\,, \eqno(4.17)$$
	  so that
	  $$ m(t)= \sum_{n=1}^\infty \int_0^t q_n(t')\, dt' = 
	  \sum_{n=1}^\infty  G_\beta \left(n^{-1/\beta} t \right)\,. \eqno(4.18)$$
	  We do not know an explicit expression for this sum if $0<\beta<1$. However, in the limiting case
	  $\beta=1$ we obtain
	  $$ m(t) = [t]=N(t)\,. \eqno(4.19)$$
	  Using (4.17) we investigate the asymptotic behaviour of $m(t)$ for $t\to \infty$.
	  We have for $s\to 0$ $\widetilde m(s)\sim 1/s^{1+\beta}$ and thus, by Tauber theory,
	  see e.g. Feller (1971) \cite{Feller_1971},
	  $$ m(t)\sim \frac{t^\beta}{\Gamma(1+\beta)}\; \hbox{for}\; t\to \infty\,. \eqno(4.20)$$
	  Remember, for the fractional Poisson process, we had found  
	 $$ m(t)= \frac{t^\beta}{\Gamma(1+\beta)}\; \hbox{for all}\; t \ge 0\,. \eqno(4.21)$$
	 {\bf Remark} 
	 A rescaled version of the discretized stable subordinator 
	 can be used for producing closely spaced precise snapshots of a true particle trajectory of a space-time 
	 fractional diffusion process, see e.g. the recent chapter by Gorenflo and Mainardi (2011)
	 \cite{Gorenflo-Mainardi_METZLER11} on parametric subordination.

\section{The diffusion limits for the fractional \\ Poisson and the Wright processes} 
In a CTRW we can,	  with positive scaling factor $h$ and $\tau$, 
replace the jumps $X$ by jumps $X_h=h\,X$,
the waiting times $T$ by waiting times $T_\tau= \tau\, T$.
This leads to the rescaled jump density 
$w_h(x)=w(x/h)/h$ and the rescaled waiting time density 
$\phi_\tau(t)= \phi(t/\tau)/\tau $
and correspondingly to the transforms
$\widehat w_h(\kappa) = \widehat w (h \kappa)$, $\widetilde \phi _\tau(s)=\widetilde \phi(\tau s) $.	 
\vsp
For the sojourn density $u_{h,\tau}(x,t)$, density in $x$ evolving in $t$, we obtain from (2.9)
in the transform domain
$$
\widehat{\widetilde{u}}_{h,\tau} (\kappa,s)
= \frac{1-\widetilde{\phi}(\tau s)}{s} \, \frac{1}{1- \widetilde{\phi}(\tau s)\,\widehat{w}(h \kappa)}
\,,
\eqno(5.1)$$
where, if $w(x)$ has support on $x\ge 0$ we can work with the Laplace transform instead of the Fourier transform
(replace the $\widehat{\phantom{x}}\,$ by $\widetilde{\phantom{x}}\,$). 
If there exists between $h$ and $\tau$ a scaling relation ${\mathcal{R}}$ (to be introduced later) 
under which $u(x,t)$ tends for $h \to 0$, $\tau \to 0$ to a meaningful limit
$v(x,t)=u_{0,0}(x,t)$, then we call the process $x=x(t)$ with this sojourn density a {\it diffusion limit}.
We find it via
$$ \widehat{\widetilde{v}} (\kappa,s)= \lim_{h, \tau \to 0 ({\mathcal{R}})}  
\widehat{\widetilde{u}}_{h,\tau} (\kappa,s)\,, \eqno(5.2)$$
and Fourier-Laplace (or Laplace-Laplace) inversion.
\vsp
Note: this diffusion limit is a limit in the weak sense (convergence in distribution of the CTRW to the diffusion 
limit). The mathematical background consists in the application of the Fourier (or Laplace) continuity 
theorem of probability theory for fixed time $t$.
\vsp
We will now find that the counting numbers of the fractional Poisson process and the Wright  process 
have the same diffusion limit, namely the inverse stable subordinator.
The two corresponding Erlang processes have the same diffusion limit, namely the stable subordinator.
For $t \to \infty$ the renewal functions have the same asymptotic behaviour, namely 
$m(t)\sim t^\beta/\Gamma(1+\beta)$.
Here,  in the case of the fractional Poisson process, we can replace the sign $\sim$ of asymptotics by the sign $=$ of equality  for all $t\ge 0$.
  \vsp
To prove these statements we need the Laplace transform of the relevant functions $\phi(t)$ and $w(x)$.
For the fractional Poisson process we have
$$ \phi(t)= \frac {d}{dt} E_\beta(-t^\beta) \, \div \, \widetilde \phi(s)= \frac{1}{1+s^\beta}\,, \q
 w(x) =\delta (x-1) \, \div \, \widetilde w(\kappa) = \exp (-\kappa)\,.$$
For the Wright process we have
 $$ \phi(t)= g_\beta(t) \,  \div \,\widetilde \phi(s) = \exp (-s^\beta)\,, \q
 w(x) =\delta (x-1) \, \div \, \widetilde w(\kappa) = \exp (-\kappa)\,.$$
 In all cases we have, for fixed $s$ and $\kappa$ 
 $$  \widetilde \phi(\tau s) \sim 1 -(\tau s)^\beta \; \hbox{as} \; \tau \to 0\,, \q
   \widetilde w(h\kappa) \sim 1 - (h\kappa)\; \hbox{as} \; h \to 0\,, \,.$$
   and straightforwardly we obtain for the sojourn densities in both cases, 
   by use of (5.1) with $p$ in place of $u$ and 
    $\widehat{\phantom{x}}\,$
    replaced by  $\widetilde{\phantom{x}}\,$
   $$
\widetilde{\widetilde{p}}_{h,\tau} (\kappa,s)
\sim \frac{\tau^\beta\, s^{\beta-1}}{\tau^\beta \,s^\beta + h\,\kappa}\,, \q
\hbox{for}\;  \tau \to 0\,, \; h \to 0\,.\eqno(5.3)$$
Using the scaling relation ${\mathcal{R}}$
$$ h = \tau^\beta\,,\eqno(5.4)$$
we obtain 
 $$
\widetilde{\widetilde{p}}_{0,0} (\kappa,s)
= \frac{ s^{\beta-1}}{s^\beta + \kappa}\,, \eqno(5.5)
  $$
  By partial Laplace inversions we get two equivalent representations
  $$ p_{0,0}(x,t) = \L_\kappa^{-1} \left\{E_\beta(-\kappa t^\beta)\right\}
 = \L_s^{-1} \left\{s^{\beta-1}\, \exp(-x s^\beta)\right\}\,,\eqno(5.6)$$
  leading to the density of the {\it inverse stable subordinator} 
   $$ p_{0,0}(x,t) = t^{-\beta}\, M_\beta(x/t^\beta) =J_t^{1-\beta} f(t,x)\,,
   \eqno(5.7)$$ 
   where $M_\beta$ and $J^{1-\beta}_t$ denote respectively the $M$-Wright function and the Riemann-Liouville 
   fractional integral introduced in Appendix A, and $f(t,x)$  
     the stable subordinator given by Eq. (4.7).
    \vsp
{\bf Remark:} In (4.7 and (5.7) the densities of the stable and the inverse stable subordinator 
are both represented via the $M$-Wright function.
\vsp
\paragraph{The diffusion limit for the Erlang process.} 
\vsp
In the Erlang process the roles of space and time, likewise of jumps and waiting times, are interchanged. 
In other words we treat $x\ge 0$ as a pseudo-time variable and $t\ge 0$ as a pseudo-space variable.
For the resulting sojourn density $q(t,x)$, we have  from interchanging in (5.1)
for $ h \to 0$ and $\tau \to 0$,
$$
\widetilde{\widetilde{q}}_{h,\tau} (s, \kappa)
= \frac{1-\widetilde{w}(h \kappa)}{k} \, \frac{1}{1- \widetilde{w}(h \kappa)\,\widetilde{\phi}(\tau s)}
\sim \frac{h}{h\kappa + (\tau s)^\beta} \; \
\,.
\eqno(5.8)$$
Again using the scaling relation ${\mathcal{R}}$ in  Eq. (5.4) we find
$$
\widetilde{\widetilde{q}}_{0,0} (s, \kappa) =\frac{1}{\kappa +s^\beta}\,,\eqno(5.5')$$
which is the Laplace-Laplace transform of the density of  stable subordinator of Section 4.
In fact, by partial Laplace inversion, 
$$\widetilde{q}_{0,0} (s, x)=\exp(-x s^\beta) =\widetilde f(s,x)\,, \eqno(5.9)$$
and it follows that 
$$ {q}_{0,0}(t,x)= f(t,x)\,, \q x\ge 0\,,\; t\ge 0\,.\eqno(5.10)$$
See (4.7) for its explicit representation as a rescaled stable density expressed via a $M$-Wright function.
\vsp
We get the same result by continualization of the discretized stable subordinator.
Replace in Eqs. (4.9), (4.10) the discrete variable $n$ by the continuous variable $x$.

\section{Conclusions}
The fractional Poisson process and the Wright process (as discretization of the stable subordinator)
along with their diffusion limits play eminent roles in theory and simulation of fractional diffusion processes.
Here we have analyzed these two processes, concretely the corresponding counting number and 
Erlang processes, the latter being  the processes inverse to the former.
Furthermore we have obtained the diffusion limits of all these processes by well-scaled refinement 
of waiting times  and jumps.    

\section*{Acknowledgements}

The authors are grateful to Professor Mathai for several invitations to visit the Centre for Mathematical Sciences in  Pala-Kerala for conferences, teaching and research. 
They luckily enjoyed there the friendly and stimulating  environment, scientifically and geographically.
The first-named author appreciates the stimulating working conditions he enjoyed during several ERASMUS visits in the Department of Physics of Bologna University.

\section*{Appendix A: Operators, transforms and special functions}\label{Appendix}
 For the reader's convenience  here we present a brief introduction to the basic notions required  
 for the presentation and analysis of the renewal processes to be treated, 
 including essentials on fractional calculus and special functions of Mittag-Leffler and Wright type. .   
\vsp
Thereby we follow our earlier papers concerning related topics, see  
\cite{Gorenflo_PALA10,
GAR_Vietnam03,%
Gorenflo-Mainardi_CISM97,%
Gorenflo-Mainardi_BAD-HONNEF08,
Gorenflo-Mainardi_EPJ-ST11,
Gorenflo-Mainardi_METZLER11,
Gorenflo-Mainardi_KILBAS12,
Gorenflo-et-al_Math-Fin01,
GorMaiViv_CSF07,
Mainardi-et-al_VIETNAM04,
Mainardi-Luchko-Pagnini_FCAA01,
Mainardi-Mura-Pagnini_IJDE10,
Scalas-et-al_PhisicaA00,
Scalas_PRE04},
and our recent monograph on Mittag-Leffler Functions and Related Topics
\cite{GKMR_BOOK14}.
 
\vsp
For more details on general aspects the interested reader may consult the  treatises, listed in order of publication time, 
by Podlubny  \cite{Podlubny_BOOK99},
Kilbas and Saigo \cite{Kilbas-Saigo_BOOK04},
Kilbas, Srivastava and Trujillo  \cite{Kilbas-et-al_BOOK06}, 
Mathai and Haubold  \cite{Mathai-Haubold_BOOK08},
Mathai, Saxena and Haubold  \cite{Mathai-Saxena-Haubold_BOOK-H-2010}, 
Mainardi  \cite{Mainardi_BOOK10}, 
Diethelm  \cite{Diethelm_BOOK10},
 Baleanu, Diethelm, Scalas and Trujillo \cite{BDST_BOOK12}, 
 Uchaikin \cite{Uchaikin_BOOK13},
Atanackovi\'c,   Pilipov\'ic,  Stankovi\'c  and Zorica \cite{Atanackovic_BOOK1-14}.
  \vsp
 \textbf{Fourier and Laplace transforms}
 \vsp
By $\RR$ ($\RR^+$, $\RR_0^+$) we mean the set of all (positive, non-negative) real numbers, and by $\CC$
the set of  complex numbers. 
It is known that the Fourier transform is applied to functions defined in $L_1(\RR\,)$ whereas
 the Laplace transform is applied to functions defined in $L_{loc}(\RR^+)$.
In our cases the arguments of the original function are the space--coordinate $x$ 
($x\in \RR\,$ or $x\in \RR_0^+$)
 and the time--coordinate $t$ ($t \in \RR_0^+$).
  We  use the symbol $\div$ for the juxtaposition of
a function  with its Fourier or Laplace transform.
A look at  the superscript $\,\widehat{\phantom{}}\,$  for the Fourier transform, 
$\,\widetilde{\phantom{}}\,$ for the Laplace transform reveals their relevant 
juxtaposition.  
We use $x$ as argument (associated to real $\kappa$) for functions Fourier transformed, and 
$x$ or $t$ as argument (associated to complex $\kappa$ or $s$, respectively) for functions Laplace transformed. 
\\
$$  f(x) \,\div\, \widehat f(\kappa) := \int_{-\infty}^{+\infty} \!\! \e^{i\kappa x}\, f(x)\, dx\,,
\q \hbox {Fourier transform}. $$
 $$ f(x) \,\div\, \widetilde f(\kappa) := \int_{0}^{\infty} \!\! \e^{-\kappa x}\,
 f(x)\, dx\,, \q \hbox {Laplace transform}.$$
 $$  f(t) \,\div\, \widetilde f(s) := \int_{0}^{\infty}\!\! \e^{-s t}\,f(t)\, dx\,,
  \q \hbox {Laplace transform}.$$
  \textbf{Convolutions}
  $$ (u*v)(x) := \int_{-\infty}^{+\infty}\!\! u(x-x')\, v(x')\, d x'\,, \q \hbox{Fourier convolution}.$$
  $$ (u*v)(t) := \int_{0}^{t} \!\! u(t-t')\, v(t')\, d t'\,, \q \hbox{Laplace convolution}.$$
  The meaning of the connective $*$
 will be clear from the context. For convolution powers we have:
 $$ u^{*0}(x) = \delta(x)\,, \; u^{*1}(x)= u(x)\,, \; u^{*(n+1)}(x)= (u^{*n}*u)(x)\,,$$
 $$ u^{*0}(t) = \delta(t)\,, \; u^{*1}(t)= u(t)\,, \; u^{*(n+1)}(t)= (u^{*n}*u)(t)\,,$$
 where $\delta$ denotes the Dirac generalized function.
 \vsp
 \textbf{Fractional integral} 
 \vsp
    The \underline{Riemann-Liouville  fractional integral} of order $\alpha>0$,
 for a sufficiently well-behaved  function $f(t)$ ($t\ge 0$), 
  is defined as  
  $$J^\alpha_t f(t)= 
  \rec{\Gamma(\alpha  )}\,
  \int_0^t (t-\tau )^{\alpha  -1}\, f(\tau )\,{\mbox{d}}\tau \,, \q \alpha>0\,,$$
  by convention as $f(t)$ for $\alpha=0$. 
 Well known are the {\it semi-group property} 
$$
J_t^\alpha  \,J_t^\beta = J_t^{\alpha  +\beta} = J_t^\beta  \,J_t^\alpha \,,
   \q \alpha \,,\;\beta  \ge 0\,,$$
and the  Laplace transform pair 
$$ J_t^\alpha  \;f(t) \div
     \frac{\widetilde f(s)}{s^\alpha }\,,\q \alpha  \ge 0\,.  $$
  \vsp
 \textbf{Fractional derivatives} 
 \vsp
  The  \underline{Riemann-Liouville  fractional derivative operator} of order $\alpha>0$,
 $D^\alpha _t $, is
 defined as the {\it left inverse operator} of the corresponding fractional integral $J^\alpha _t$.
  Limiting ourselves to fractional derivatives of order $\alpha \in (0,1)$ we have, 
 for a sufficiently well-behaved  function $f(t)$ ($t\ge 0$),   
$$D^\alpha_t  \,f(t) \!:= D_t^1\, J_t^{1-\alpha}\, f(t) =
  {\ds \rec{\Gamma(1-\alpha )}}
  {\ds \frac{{\mbox{d}}}{{\mbox{d}}t}
  \int_0^t 
    \frac{f(\tau)}{(t-\tau )^{\alpha }}\, {\mbox{d}}\tau}\,, \; 
0  <\alpha  < 1 \,, $$ 
while the corresponding \underline{Caputo derivative} is
 $$
	 \begin{array}{ll}
	  {\ds _*D^\alpha_t  \,f(t)} & \!\! := \! J_t^{1-\alpha}\, D_t^1\,f(t) \! = \!
	   {\ds \rec{\Gamma(1-\alpha )}} 
  {\ds  \int_0^t  \frac{f^{(1)}(\tau)}{(t-\tau )^{\alpha}} \,{\mbox{d}}\tau }\\
  & \!\! = \!\!
  {\ds D^\alpha_t f(t) -  f(0^+) \frac{t^{-\alpha }} {\Gamma(1-\alpha )} =\
 D^\alpha_t  \left[ f(t) - f(0^+) \right]} \,.
  \end{array} 
  $$ 
  Both derivatives yield the ordinary first derivative as $\alpha \to 1^-$ but for $\alpha \to 0^+$
  we have
  $$ D^0_t\,f(t)= f(t)\,, \quad  \,_*D^0_t  \,f(t) = f(t) -f(0^+)\,.$$
 We point out the major utility
of the Caputo fractional derivative  
in treating initial-value problems with Laplace transform. We have
$$  \L[\,_*D^\alpha_t\, f(t); s] 
= s^\alpha \widetilde f(s)- s^{\alpha-1}\, f(0^+)\,, \q 0<\alpha\le 1\,. $$
In contrast the Laplace transform of the Riemann-Liouville fractional derivative 
needs the limit at zero of a fractional integral of the function $f(t)$.
\vsp
Note that both types of fractional derivative may exhibit singular behaviour at the origin $t=0^+$. 
   \vsp
 \textbf{Mittag-Leffler and Wright functions}
 \vsp
 \underline{The Mittag-Leffler function of parameter $\alpha$} is defined as
 $$
 E_\alpha (z) := \sum_{n=0}^\infty \frac{z^n}{\Gamma (\alpha n+1)}
\,,\q \alpha > 0\,, \q z\in \CC\,.$$
It is  entire of order $1/\alpha$.
Let us note the trivial cases 
$$ \left\{
\begin{array}{ll}
E_1(\pm z)= \exp (\pm z)\,, \\
 E_2\left(+z^2\right) =  \cosh \, (z)\,, \;
   E_2\left(-z^2\right) =  \cos \, (z)\,.
   \end{array}
   \right.
    $$
Without  changing the order $1/\alpha$
the Mittag-Leffler function can be  generalized by introducing  an additional (arbitrary)
parameter $\beta$.
 \vsp
\underline{The Mittag-Leffler function of parameters $\alpha, \beta$} is defined as
 $$
 E_{\alpha,\beta} (z) := \sum_{n=0}^\infty \frac{z^n}{\Gamma (\alpha n+\beta)}
\,,\q \alpha > 0\,, \q \beta, \, z\in \CC\,.$$
\underline{Laplace transforms of Mittag-Leffler functions}
\vsp
   For our purposes we  need, with $0<\nu \le 1$ and $t\ge 0$,
   the Laplace transform pairs 
  $$ \left\{
  \begin{array}{ll}
\Psi(t) \!=\! E_\beta(-t^\nu) \, \div\, {\ds \widetilde \Psi(s)= \frac{s^{\nu-1}}{1+ s^\nu}}\,,
\\
\phi(t)  \!= \!- {\ds \frac{d}{dt}} E_\nu(-t^\nu) = \!t^{\nu-1}\, E_{\nu,\nu}(-t^\nu)
	\, \div \, {\ds \widetilde \phi(s) \!= \!\frac{1}{1+ s^\nu}}\,.
\end{array}
	\right. $$
	Of high relevance is  
the algebraic decay of $\Psi(t)$ and $\phi (t)$ 
as $t\to \infty$:  
$$ \left\{
\begin{array}{ll}
{\ds \Psi(t)\sim    \frac{\sin (\nu  \pi)}{\pi}\,\frac{\Gamma(\nu )}{t^\nu }}\,,
\\ 
 {\ds \phi  (t) 
   \sim   \frac{\sin (\nu  \pi)}{\pi}\,\frac{\Gamma(\nu +1)}{t^{\nu +1}}}\,,
 \end{array}
 \right.
    \;  t \to +\infty\,.
	$$
Furthermore  $\Psi(t) = E_\nu(-t^\nu)$ is the solution of the fractional  relaxation equation
with the Caputo derivative
$$ 
        _*D^\nu _t u(t)  = - u(t)\,, \q t\ge 0\,, \quad u(0^+) =1\, $$
whereas $\phi(t) = -{\ds\frac{d}{dt} E_\nu(-t^\nu)} $ is the solution of the fractional  relaxation equation
with  the Riemann-Liouville derivative
$$D^\nu_t u(t)  = - u(t)\,, \q t\ge 0\,, \quad \lim_{t\to 0^+} J_t^{1-\nu}u(t) =1\,. $$
We refer to our recent monograph \cite{GKMR_BOOK14} and to our papers
\cite{Gorenflo_PALA10,GAR_Vietnam03,%
Gorenflo-Mainardi_BAD-HONNEF08, Gorenflo-Mainardi_EPJ-ST11, %
Gorenflo-Mainardi_METZLER11, Gorenflo-Mainardi_KILBAS12}
for the relevance of Mittag-Leffler functions in
theory of continuous time random walk and space-time fractional diffusion
and in power law asymptotics.
Particularly worth to be mentioned  is the pioneering paper by Hilfer and Anton \cite{Hilfer-Anton_PRE95}.  They show that for transforming a general evolution equation for continuous time random walk into the time fractional version of the Kolmogorov-Feller equation a waiting time law expressible via a Mittag-Leffler type function is required.
\vsp
\underline{The Wright function} is defined as
 $$ W_{\lambda ,\mu }(z ) :=
   \sum_{n=0}^{\infty}\frac{z^n}{n!\, \Gamma(\lambda  n + \mu )}\,,
 \q \lambda  >-1\,, \,\q \mu \in \CC\,, \q z \in \CC\,.$$
 We distinguish the Wright functions of 
{\it first kind} ($\lambda \ge 0$)
and {\it second kind} ($-1<\lambda<0$).
The case $\lambda =0$ is trivial since
$  W_{0, \mu }(z) = { \e^{\, z}/ \Gamma(\mu )}\,.$
The Wright function is entire of order
$1/(1+\lambda)$ hence
of exponential type only if $\lambda \ge 0$.
\vsp
\underline{Laplace transforms of the Wright functions}
\vsp
For the Wright function of the first kind, being entire of exponential type,
the Laplace transform can be obtained  by transforming the power series term by term:
$$  W_{\lambda, \mu}(t) \,\div \, 
    \rec{s}\, E_{\lambda ,\mu }\left( \rec{s}\right) \,, \q \lambda\ge 0\,,.
 $$
 For the Wright function of the second kind, denoting $\nu= |\lambda| \in (0,1) $ we have
   with $\mu >0$ for simplicity, we have 
$$     W_{-\nu ,\mu }(-t) \,\div\,
   E_{\nu  , \mu +\nu  }(-s)\,, \q 0<\nu  <1\,.$$
We note the minus sign in the argument in order to ensure the the existence of 
the Laplace transform thanks to the Wright asymptotic formula valid in a certain sector
symmetric to and including   the negative real axis.
\vsp
\underline{Stretched Exponentials as Laplace transforms of Wright functions}
We outline the following Laplace transform pairs related to the stretched exponentials in the
transform domain, useful for our purposes,
$$ \rec{t}\, W_{-\nu,0}\left(-\rec {t^\nu} \right) \,\div\, \e ^{-s^\nu}\,,$$ 
 $$  W_{-\nu,1-\nu}\left(-\rec {t^\nu} \right) \,\div\,, \frac{\e ^{-s^\nu}}{s^{1-\nu}}\,,$$ 
 $$  W_{-\nu,1}\left(-\rec {t^\nu} \right) \,\div\, \frac{\e ^{-s^\nu}}{s}\,.$$
 For $\nu=1/2$ we have the three sister functions related to the diffusion equation available 
 in most Laplace transform handbooks
 $$  \rec{2\sqrt{\pi}}\, t^{-3/2}\, \e^{-1/(4t)}\, \div \, \e ^{-s^{1/2}}\,,$$
 $$  \rec{\sqrt{\pi}}\, t^{-1/2}\, \e^{-1/(4t)}\, \div \, \frac{\e ^{-s^{1/2}}}{s^{1/2}}\,,$$
 $$ \hbox{erfc}\left(\frac{1}{2t^{1/2}}\right)\,\div \, \frac{\e ^{-s^{1/2}}}{s}\,.$$
 Among the Wright functions of the second kind a fundamental role in fractional diffusion equations  
is played by the so called $M$-Wright function,
see e.g. \cite{Mainardi_BOOK10,Mainardi-Luchko-Pagnini_FCAA01,Mainardi-Mura-Pagnini_IJDE10}.
\vsp
\underline{The $M$-Wright function} is defined as
$$ \! M_\nu (z) \!:=\!  W _{-\nu , 1-\nu }(-z) \!=\!
 {\sum_{n=0}^{\infty}
 \frac{(-z)^n }{  n! \Gamma[-\nu n + (1-\nu )]} } 
\!=\! {\rec{\pi} \sum_{n=1}^{\infty}\frac{(-z)^{n-1} }{  (n-1)!}
  \Gamma(\nu n)  \sin (\pi\nu n)}  \,,
  $$
 with $z\in \CC$ and $0<\nu<1$. Special cases are 
 $$  M_{1/2}(z) \!=\!
  \rec{\sqrt{\pi}}\, \exp \left(-{\,z^2/ 4}\right),\;
  M_{1/3}(z)  \!=\! 3^{2/3}  {\rm Ai} \left( {z/ 3^{1/3}}\right)\,.
 $$
 where $Ai$ denotes the Airy function, see e.g.
 \cite{AS_65}.
 
 \underline{The asymptotic representation of the $M$-Wright function}	
\vsp
	Choosing as a variable $t/\nu $ rather than $t$, the computation of the 
asymptotic representation as $t\to \infty$ by the saddle-point approximation
 yields:
$$  M_\nu (t/\nu ) \sim
   a(\nu )\, {\ds t^{(\nu -1/2)/(1-\nu)}}
   \exp \left[- b(\nu)\,{\ds t^{1/(1-\nu)}}\right]\,,
$$
where $$ a(\nu) = \rec{\sqrt{2\pi\,(1-\nu)}}    >0 \,,  \q
  b(\nu) = \frac{1-\nu }{  \nu }    >0 \,.$$
 \vsp
 \underline{Mittag-Leffler function as Laplace transforms of $M$-Wright function}
    $$  M_\nu (t) \,\div\,  E_\nu (-s)\,, \q 0<\nu<1\,, \q t\ge 0\,, \q s \ge 0\,.  $$
\underline{Stretched Exponentials as Laplace transforms of $M$-Wright functions}
   $$ \frac{\nu }{  t^{\nu +1}}\,  M_\nu \left( 1/{t^\nu } \right)\,\div\,
    \e^{\ds \,-s^\nu}\,, \;  0<\nu <1\,, \q t\ge 0\,, \q s \ge 0\,. $$
$$\frac{1}{  t^{\nu}}\,  M_\nu \left( 1/{t^\nu } \right)\,\div\,
    \frac{\e^{\ds\, -s^\nu}}{s^{1-\nu}}\,, \;  0<\nu <1\,,\q t\ge 0\,, \q s \ge 0\,.  $$
	Note that $\exp (-s^\nu)$ is the Laplace transform of the extremal (unilateral) stable density
	$L_\nu^{-\nu}(t)$, which vanishes for $t< 0$, 
	so that, introducing the Riemann-Liouville fractional integral, we have
	$$\frac{1}{  t^{\nu}}\,  M_\nu \left( 1/{t^\nu } \right)
	= J_t^{1-\nu} \left\{L_\nu^{-\nu}(t)\right\} 
	= J_t^{1-\nu} \left\{\frac{\nu }{  t^{\nu +1}}\,  M_\nu \left( 1/{t^\nu } \right)\right\}\,.$$
	
    	\section*{Appendix B: Collection of results}\label{Appendix B}
\vsp	
{\bf General renewal process}
\vsp
Waiting time density: $\phi(t)$; 
Survival function: $\Psi(t)= {\ds \int_t^\infty \phi(t')\, dt'}$
\vsp
(a) \underline{The counting number process} $x=N(t)$ has probability density function 
(density in $x\ge 0$ and evolving in $t\ge 0$):
$$ p(x,t) = \sum_{n=0}^\infty p_n(t)\, \delta(x-n)\,,$$
and counting probabilities
$$ p_n(t)= \left(\Psi \,*\, \phi^{*n}\right) (t)\,.$$	
	(b) \underline{The Erlang process} $t=t(n)$, inverse to the counting process has probability density function
(density in $t$, evolving in $n=0,1,2, \dots$)	
$$q_n(t) = q(t,n) = \phi^{*n} (t)\,,$$ 
with 
$$q_n(t) = \frac{d}{dt} Q_n(t)\,, \q Q_n(t) = \sum_{k=n}^\infty p_k(t)\,,$$
where $q_n(t), \, Q_n(t)$ are the Erlang densities and  probability distribution functions, respectively.
Note that $p_n(t)= \left(\Psi \,*\, q_n\right) (t)$.
\vsp
{\bf Special cases}
\vsp 
($\alpha$) \underline{The fractional Poisson process}
$$\phi(t) =  - {\ds \frac{d}{dt} E_\beta(-t^\beta)}
\, \div\, \widetilde \phi(s)= \frac{1}{1+s^\beta}\,,$$
$$ p_n(t)= \frac{t^{n\beta}}{n!}\, E_\beta^{(n)}\left(-t^\beta\right)\,.$$
The Erlang densities are
$$ q_n(t)= \beta \frac{t^{n\beta-1}}{(n-1)!}\, E_\beta^{(n)}\left(-t^\beta\right)\,.$$
\vsp
($\beta$) \underline{The Wright process}
$$\phi(t) =  g_\beta(t)
\, \div\, \widetilde g_\beta(s)= \exp (-s^\beta)\,,$$
$$ p_n(t)= G_\beta\left(n^{-1/\beta}t\right) - G_\beta\left((n+1)^{-1/\beta}t\right)\,.$$
The Erlang densities are
$$ q_n(t)= n^{-1/\beta}\, g_\beta\left(n^{-1/\beta}\,t\right)\,.$$
\vsp

\end{document}